\documentclass[a4paper,10pt]{siamltex}
\usepackage[dvips]{epsfig,graphicx}
\usepackage{verbatim}
\usepackage{array}
\usepackage{color}
\usepackage{colortbl}
\usepackage{fancybox}
\usepackage{longtable}
\usepackage{epsfig}
\usepackage{amsmath,amsfonts}
\usepackage{amssymb}
\usepackage{latexsym}
\usepackage{graphicx}
\usepackage{rotating}
\usepackage{longtable,dcolumn,multirow}

\newcommand{\norm}[1]{\Vert {#1} \Vert}

\newcommand{\neuc}[1]{\left \Vert {#1} \right \Vert_2}

\newcommand{\nfro}[1]{\left \Vert {#1} \right \Vert_F}

\newcommand{\da}{\Delta A}
\newcommand{\db}{\Delta b}
\newcommand{\dab}{(\Delta A,\Delta b)}

\newcommand{\bigo}{\mathcal{O} }

\newcommand{\C}{\def\C{\mbox{I\hspace{-.47em}C}}}
\newtheorem{theo}{Theorem}

\newtheorem{Remark}{Remark}
\newtheorem{Algo}{Algorithm}

\title{Computing the Conditioning of the Components of a Linear Least Squares Solution}
\author{
Marc Baboulin\footnotemark[1]\thanks{CERFACS, and University of Tennessee, 
   Email~: {\tt baboulin@eecs.utk.edu}}
\and Jack Dongarra\footnotemark[2]\thanks{University of Tennessee, Oak Ridge National 
   Laboratory, and University of Manchester, Email~: {\tt dongarra@eecs.utk.edu}}
\and Serge Gratton\footnotemark[3]\thanks{Centre National d'Etudes Spatiales and CERFACS,
   Email~: {\tt gratton@cnes.fr}}
\and Julien Langou\footnotemark[4]\thanks{
   University of Colorado at Denver and Health Sciences Center, 
        Email~: {\tt julien.langou@cudenver.edu}}
}
\begin{document}
\thispagestyle{empty}
\maketitle
\begin{abstract}

In this paper, we address the accuracy of the results for the overdetermined
full rank linear least squares problem.  We recall theoretical results obtained
in~\cite{ABG.07} on conditioning of the least squares solution and the
components of the solution when the matrix perturbations are measured in
Frobenius or spectral norms.
Then we define computable estimates for these condition numbers and we interpret 
them in terms of statistical quantities.  In particular, we show that, in the
classical linear statistical model, the ratio of the variance of one component
of the solution by the variance of the right-hand side is exactly the condition
number of this solution component when perturbations on the right-hand
side are considered. We also provide fragment codes using LAPACK~\cite{LAPACK} routines 
to compute the variance-covariance matrix 
and the least squares conditioning 
and we give the corresponding computational cost. 
Finally we present a small historical numerical example that was used by 
Laplace~\cite{LAPLACE} for computing the mass of Jupiter 
and experiments from the space industry with real physical data.\\


\end{abstract}

\renewcommand{\thefootnote}{\fnsymbol{footnote}}

\section{Introduction}

We consider the linear least squares problem (LLSP)
$\min_{x \in \mathbb{R}^n}\|Ax-b\|_2$, where $b \in \mathbb{R}^m$
and $A \in \mathbb{R}^{m\times n}$ is a matrix of full column rank $n$.\\
Our concern comes from the following observation: in many parameter estimation
problems, there may be random errors in the observation vector $b$ due to
instrumental measurements as well as roundoff errors in the algorithms. 
The matrix $A$ may be subject to errors in its computation 
(approximation and/or roundoff errors). 
In such cases, while the condition number of the matrix 
$A$ provides some information
about the sensitivity of the LLSP to perturbations, a single global
conditioning quantity is often not relevant enough since we may have
significant disparity between the errors in the solution components.
We refer to the last section of the manuscript for
illustrative examples.\\
There are several results for analyzing
the accuracy of the LLSP by components.
For linear systems $Ax=b$ and for LLSP,
~\cite{CHA.IPS.95} defines so called componentwise condition numbers that correspond 
to amplification factors of the relative errors in solution components due to
perturbations in data $A$ or $b$ and explains how to estimate them.
For LLSP,~\cite{KEN.LAU.98} proposes to estimate componentwise condition
numbers by a statistical method. More recently,~\cite{ABG.07} developed
theoretical results on conditioning of linear functionals of
LLSP solutions.\\
%
%
%
The main objective of our paper is to provide computable quantities of these
theoretical values in order to assess the accuracy of an LLSP solution or some
of its components.  To achieve this goal, traditional tools for the numerical
linear algebra practitioner are condition numbers or backward errors whereas
the statistician usually refers to variance or covariance.  Our purpose here is
to show that these mathematical quantities coming either from numerical
analysis or statistics are closely related.  In particular, we will show in
Equation~(\ref{eq:1}) that, in the classical linear statistical model, the
ratio of the variance of one component of the solution by the variance of the
right-hand side is exactly the condition number of this component 
when perturbations on the right-hand side are considered.  In that
sense, we attempt to clarify, similarly to~\cite{HIGH.STEW.87}, the analogy
between quantities handled by the linear algebra and the statistical approaches
in linear least squares.  Then we define computable estimates for these
quantities and explain how they can be computed using the standard libraries
LAPACK or ScaLAPACK.
%

This paper is organized as follows.
In Section~\ref{sec:theoback}, we recall and exploit some results of practical interest
coming from~\cite{ABG.07}. We also define the condition
numbers of an LLSP solution or one component of it.
In Section~\ref{sec:link}, we recall some definitions and results related to
the linear statistical model for LLSP, and we interpret
the condition numbers in terms of statistical quantities.
In Section~\ref{sec:lapack} we provide practical formulas and
FORTRAN code fragments for computing
the variance-covariance matrix and LLSP condition numbers using LAPACK.
In Section~\ref{sec:numerics}, we propose two numerical examples
that show the relevance of the proposed quantities and their practical
computation. The first test case is a historical example from Laplace
and the second example is related to gravity field computations.
Finally some concluding remarks are given in Section~\ref{sec:concl}.\\

Throughout this paper we will use the following notations.
We use the Frobenius norm $\nfro{.}$ 
and the spectral norm $\neuc{.}$ on matrices 
and the usual Euclidean norm $\neuc{.}$ on vectors. 
$A^{\dagger}$ denotes the Moore-Penrose pseudo inverse of $A$,
the matrix $I$ is the identity matrix and $e_i$ is the $i$-th canonical vector.
\section{Theoretical background for linear least squares conditioning}\label{sec:theoback}
Following the notations in~\cite{ABG.07}, we consider the function
\begin{equation}
\label{funct}
\begin{array}{c c c c}
g\ : &
\mathbb{R} ^{m \times n} \times \mathbb{R} ^m &  \longrightarrow &  \mathbb{R} ^k
\\
 &  A,b & \longmapsto &  g(A,b)=L^{T}x(A,b)=L^T(A^TA)^{-1}A^Tb,\\
 \end{array}
\end{equation}
where $L$ is an $n \times k$ matrix, with $k \leq n$. 
Since $A$ has full rank $n$, $g$ is continuously F-differentiable in a
neighbourhood of $(A,b)$ and we
denote by $g'$ its F-derivative.\\
Let $\alpha$ and $\beta$ be two positive real numbers.
In the present paper we consider 
the Euclidean norm for the solution space $\mathbb{R}^k$. 
For the data space $\mathbb{R}^{m\times n}\times \mathbb{R}^m$, 
we use the product norms defined by
$$\norm{(A,b)}_{\rm{F~or~2}}=
\sqrt{\alpha^2\norm{A}_{\rm{F~or~2}}^2+\beta^2\neuc{b}^2},~\alpha,\beta>0.$$
%
%
%
Following~\cite{GEURTS}, the absolute condition number of $g$ at the point
$(A,b)$ using the product norm defined above is given by:
$$\kappa_{g,{\rm{F~or~2}}}(A,b)
=\max_{\dab}
\frac{\neuc{g'(A,b).\dab}}{\norm{\dab}_{\rm{F~or~2}}}.$$
The corresponding relative condition number of $g$ at $(A,b)$ is expressed by
$$\kappa_{g,{\rm{F~or~2}}}^{(rel)}(A,b)=
\frac{\kappa_{g,F}(A,b)~\norm{(A,b)}_{\rm{F~or~2}}}{\neuc{g(A,b)}}.$$
To address the special cases where only $A$ (resp. $b$) is perturbed, 
we also define the quantities 
$\kappa_{g,{\rm{F~or~2}}}(A)
=\max_{\da}
\frac{\neuc{\frac{\partial g}{\partial A}(A,b).\da}}{\norm{\da}_{\rm{F~or~2}}}$
(resp. $\kappa_{g,2}(b)
=\max_{\db}
\frac{\neuc{\frac{\partial g}{\partial b}(A,b).\db}}{\norm{\db}_2}$).
\begin{Remark}
{\em
The product norm for the data space is very flexible; 
the coefficients $\alpha$ and $\beta$ allow us to monitor
the perturbations on $A$ and $b$.
For instance, large values of $\alpha$ (resp. $\beta$ ) enable us to obtain
condition number problems where mainly $b$ (resp. $A$) are perturbed.
In particular, we will address the special cases where 
only $b$ (resp. $A$) is perturbed by choosing the $\alpha$ and $\beta$
parameters as  
$\alpha=+\infty~{\rm and}~\beta=1$
(resp. $\alpha=1~{\rm and}~\beta=+\infty$) since we have
$$\lim_{\alpha \rightarrow +\infty}\kappa_{g,\rm{F~or~2}}(A,b)
=\frac{1}{\beta}\kappa_{g,2}(b)
~{\rm and}~
\lim_{\beta \rightarrow +\infty}\kappa_{g,\rm{F~or~2}}(A,b)
=\frac{1}{\alpha}\kappa_{g,\rm{F~or~2}}(A).
$$
This can be justified as follows:
\begin{eqnarray*}
\kappa_{g,\rm{F~or~2}}(A,b) & = &
\max_{\dab} 
\frac{\neuc{\frac{\partial g}{\partial A}(A,b).\da
+\frac{\partial g}{\partial b}(A,b).\db}}
{\sqrt{\alpha^2\norm{\da}_{\rm{F~or~2}}^2+\beta^2\neuc{\db}^2}}\\
 & = & 
\max_{\dab} 
\frac{\neuc{\frac{\partial g}{\partial A}(A,b).\frac{\da}{\alpha}
+\frac{\partial g}{\partial b}(A,b).\frac{\db}{\beta}}}
{\sqrt{\norm{\da}_{\rm{F~or~2}}^2+\neuc{\db}^2}}.\\
\end{eqnarray*}
The above expression represents the operator norm of a linear functional 
depending continuously on $\alpha$, and then we get 
$$\lim_{\alpha \rightarrow +\infty}\kappa_{g,\rm{F~or~2}}(A,b)
=\max_{\dab} 
\frac{\neuc{\frac{\partial g}{\partial b}(A,b).\frac{\db}{\beta}}}
{\sqrt{\norm{\da}_{\rm{F~or~2}}^2+\neuc{\db}^2}}
=\max_{\db} 
\frac{\neuc{\frac{\partial g}{\partial b}(A,b).\frac{\db}{\beta}}}
{\neuc{\db}}
=\frac{1}{\beta}\kappa_{g,2}(b).
$$
The proof is the same for the case where $\beta=+\infty$.\\
}
\end{Remark}\\
The condition numbers related to $L^Tx(A,b)$ are referred to as {\bf partial condition
numbers} (PCN) of the LLSP with respect to the linear operator $L$
in~\cite{ABG.07}.\\
%
%
We are interested in computing the PCN for two special cases.
The first case is when $L$ is
the identity matrix (conditioning of the solution) and the second case is when $L$ is a canonical
vector $e_i$ (conditioning of a solution component).  We can extract
from~\cite{ABG.07} two theorems that can lead to computable quantities in these
two special cases.\\
\begin{theo} \label{theobound}
In the general case where~$(L \in \mathbb{R}^{n\times k})$, the absolute
condition numbers of $g(A,b)=L^Tx(A,b)$ in the Frobenius and spectral norms can
be respectively bounded as follows
$$\frac{1}{\sqrt{3}}f(A,b) \leq \kappa_{g,F}(A,b) \leq f(A,b)$$
$$\frac{1}{\sqrt{3}}f(A,b) \leq \kappa_{g,2}(A,b) \leq \sqrt{2} f(A,b)$$
where
\begin{equation}\label{eq:equationforf(A,b)}
f(A,b)=
\left(\neuc{L^T(A^TA)^{-1}}^2 \frac{\neuc{r}^2}{\alpha^2}
+\neuc{L^T A^{\dagger}}^2 (\frac{\neuc{x}^2}{\alpha^2}
+\frac{1}{\beta^2})\right)^{\frac{1}{2}}.
\end{equation}
\end{theo}
\begin{theo}\label{corocond}
In the two particular cases:
\begin{enumerate}
\item $L$ is a vector ($L \in \mathbb{R}^n$), or
\item $L$ is the $n$-by-$n$ identity matrix ($L=I$)
\end{enumerate}
the absolute condition number of $g(A,b)=L^Tx(A,b)$ in the Frobenius norm is
given by the formula:
$$
\kappa_{g,F}(A,b)=\left(\neuc{L^T(A^TA)^{-1}}^2 \frac{\neuc{r}^2}{\alpha^2}
+\neuc{L^T A^{\dagger}}^2 (\frac{\neuc{x}^2}{\alpha^2}
+\frac{1}{\beta^2})\right)^{\frac{1}{2}}.
$$
\end{theo}
Theorem~\ref{corocond} provides the exact value for the condition number in the
Frobenius norm for our two cases of interest ($L=e_i$ and $L=I$).  From
Theorem~\ref{theobound}, we observe that
\begin{equation}\label{eq:condfrob_or_condspec}
\frac{1}{\sqrt{3}} \kappa_{g,F}(A,b)
\leq 
\kappa_{g,2}(A,b)
\leq
\sqrt{6}\kappa_{g,F}(A,b).
\end{equation}
which states that the partial condition number in spectral norm is of the same
order of magnitude as the one in Frobenius norm. In the remainder of the
paper, the focus is given to the partial condition number in Frobenius norm only.\\
For the case $L=I$, the result of Theorem~\ref{corocond} is similar
to~\cite{GR.96} and~\cite[p. 92]{GEURTS}.  The upper bound for
$\kappa_{2,F}(A,b)$ that can be derived from
Equation~(\ref{eq:condfrob_or_condspec}) is also the one obtained by~\cite{GEURTS}
when we consider pertubations in $A$.\\
%
%
%
Let us denote by $\kappa_i(A,b)$ the condition number related to the component
$x_i$ in Frobenius norm (i.e $\kappa_i(A,b)=\kappa_{g,F}(A,b)$ where
$g(A,b)=e_i^Tx(A,b)=x_i(A,b)$).  Then replacing $L$ by $e_i$ in
Theorem~\ref{corocond} provides us with an exact expression for computing
$\kappa_i(A,b)$, this gives
\begin{equation}\label{eq:componentwise_formula}
\kappa_i(A,b)=\left(\neuc{e_i^T(A^TA)^{-1}}^2 \frac{\neuc{r}^2}{\alpha^2}
+\neuc{e_i^T A^{\dagger}}^2 (\frac{\neuc{x}^2}{\alpha^2}
+\frac{1}{\beta^2})\right)^{\frac{1}{2}}.
\end{equation}
$\kappa_i(A,b)$ will be referred to as {\bf the condition number of the
solution component $x_i$}.\\
%
%
Let us denote by $\kappa_{LS}(A,b)$ the condition number related to the
solution $x$ in Frobenius norm (i.e $\kappa_{LS}(A,b)=\kappa_{g,F}(A,b)$ where
$g(A,b)=x(A,b)$).  Then Theorem~\ref{corocond} provides us with an exact
expression for computing $\kappa_{LS}(A,b)$, that is
\begin{equation}\label{eq:normwise_formula}
\kappa_{LS}(A,b)=
   \neuc{(A^TA)^{-1}}^{1/2}
   \left(
      \frac{
         \neuc{(A^TA)^{-1}} \neuc{r}^2 + \neuc{x}^2
      }
      {
         \alpha^2
      }
      + \frac{1}{\beta^2}
   \right)^{\frac{1}{2}}.
\end{equation}
where we have used the fact that $\neuc{ (A^T A )^{-1}} = \neuc{ A^\dagger }^2$.\\
$\kappa_{LS}(A,b)$ will be referred to as {\bf the condition number of the
least squares solution}.\\
Note that~\cite{IRLLS.07} defines condition numbers for both $x$ and $r$ 
in order to derive error bounds for $x$ and $r$ but uses infinity-norm 
to measure perturbations.\\
In this paper, we will also be interested in the special case where only
$b$ is perturbed ($\alpha=+\infty$ and $\beta=1$). In this case, we
will call $\kappa_i(b)$ the condition number of the solution component $x_i$,
and $\kappa_{LS}(b)$ the condition number of the least squares solution.
When we restrict the perturbations to be on $b$,
Equation~(\ref{eq:componentwise_formula}) simplifies to
\begin{equation}\label{eq:justb_componentwise}
\kappa_i(b)=\neuc{e_i^T A^{\dagger}},
\end{equation}
and Equation~(\ref{eq:normwise_formula}) simplifies to
\begin{equation}\label{eq:justb_normwise}
\kappa_{LS}(b)= \neuc{A^{\dagger}}.
\end{equation}
This latter formula is standard and is in accordance with
~\cite[p. 29]{BJORCK}.
\section{Condition numbers and statistical quantities}\label{sec:link}
\subsection{Background for the linear statistical model}\label{sec:cov}
We consider here the classical linear statistical model
$$
b = Ax +\epsilon,
~A \in \mathbb{R}^{m\times n},
~b \in \mathbb{R}^m,
{\rm rank}(A)=n,
$$
where $\epsilon$ is a vector of random errors having expected value
$E(\epsilon)=0$ and variance-covariance $V(\epsilon)=\sigma_b^2I$.
In statistical language, the matrix $A$ is referred to as the regression matrix
and the unknown vector $x$ is called the vector of regression coefficients.\\
Following the Gauss-Markov theorem~\cite{ZELEN},
the least squares estimates $\hat{x}$ is the linear unbiased estimator
of $x$ satisfying
$$\|A\hat{x}-b\|_2=\min_{x \in \mathbb{R}^n}\|Ax-b\|_2,$$
with minimum variance-covariance equal to
\begin{equation}\label{eq;variance-covariance}
   C=\sigma_b^2 (A^TA)^{-1}.
\end{equation}
Moreover $\frac{1}{m-n}\neuc{b-A\hat{x}}^2$ is an unbiased estimate
of $\sigma_b^2$. This quantity is sometimes called the mean squared error (MSE).\\
The diagonal elements $c_{ii}$ of $C$ give the variance of each component
$\hat{x}_i$ of the solution. The off-diagonal elements $c_{ij},~i \neq j$
give the covariance between $\hat{x}_i$ and $\hat{x}_j$.\\
We define $\sigma_{\hat x_i}$ as the standard deviation of the solution component 
$\hat x_i$ and we have
\begin{equation}\label{eq:654}
   \sigma_{\hat x_i} = \sqrt{c_{ii}}.
\end{equation}
In the next section, we will prove that the condition numbers $\kappa_{i}(A,b)$ and
$\kappa_{LS}(A,b)$ can be related to the statistical 
quantities $\sigma_{\hat x_i}$ and $\sigma_b$.
\subsection{Perturbation on $b$ only}
Using Equation~(\ref{eq;variance-covariance}), the variance $c_{ii}$ of the
solution component $\hat x_i$ can be expressed as
$$
   c_{ii} = e_i^T C e_i = \sigma_b^2 e_i^T ( A^{T} A )^{-1} e_i.
$$
We note that $ ( A^TA ) ^{-1} =  A^{\dagger}A^{\dagger T}$ so that
$$
    c_{ii}
  = \sigma_b^2 e_i^T ( A^{\dagger}A^{\dagger T} )  e_i
  = \sigma_b^2 \neuc{e_i^T A^{\dagger}}^2.
$$
Using Equation~(\ref{eq:654}), we get
$$
   \sigma_{\hat x_i} = \sqrt{c_{ii}} = \sigma_b \neuc{e_i^T A^{\dagger}}.
$$   
Finally from Equation~(\ref{eq:justb_componentwise}), we get 
\begin{equation}\label{eq:1}
   \sigma_{\hat x_i} = \sigma_b \kappa_i(b).
\end{equation}
Equation~(\ref{eq:1}) shows that the condition number $\kappa_i(b)$
relates linearly the standard deviation of $\sigma_b$ with the 
standard deviation of $\sigma_{\hat x_i}$.\\
Now if we consider the constant vector $\ell$ of size $n$, we have~(see \cite{ZELEN})
$$ {\rm variance}( \ell^T \hat x) = \ell^T C \ell. $$
Since $C$ is symmetric, we can write
$$
\max_{\| \ell \|_2=1} {\rm variance}(\ell^T \hat x) =  \neuc{C}. 
$$
Using the fact that $\| C \|_2 =\sigma_b^2 \neuc{ (A^TA)^{-1}}=\sigma_b^2 \neuc{ A^{\dagger} }^2 $, and
Equation~(\ref{eq:justb_normwise}), we get
$$
\max_{\| \ell \|_2=1} {\rm variance}(\ell^T \hat x) =  \sigma_b^2 \kappa_{LS}(b)^2
$$
or, if we call $\sigma( \ell^T \hat x ) $ the standard deviation of $\ell^T \hat x$,
$$
\max_{\| \ell \|_2=1} \sigma(\ell^T \hat x) =  \sigma_b \kappa_{LS}(b).
$$
Note that
$\sigma_b = \max_{\| \ell \|_2=1} \sigma(\ell^T \epsilon) $ since 
 $V(\epsilon)=\sigma_b^2I$.
\begin{Remark}
{\em
Matlab proposes a routine LSCOV that computes the quantities $\sqrt{c_{ii}}$ 
in a vector STDX and the mean squared error MSE using the syntax [X,STDX,MSE] = LSCOV(A,B).\\
Then the condition numbers $\kappa_i(b)$  
can be computed by the matlab expression
STDX/sqrt(MSE).
}
\end{Remark}

\subsection{Perturbation on $A$ and $b$}
We now provide the expression of the condition number provided 
in Equation~(\ref{eq:componentwise_formula}) and in
Equation~(\ref{eq:normwise_formula}) in term of statistical quantities.\\
Observing the following relations
$$ 
C_i = \sigma_b^2 e_i^T(A^TA)^{-1} \quad {\rm and} \quad
c_{ii} = \sigma_b^2 \neuc{e_i^T A^{\dagger}}^2,
$$
where $C_i$ is the $i$-th column of the variance-covariance matrix,
the condition number of $x_i$ given in Formula~(\ref{eq:componentwise_formula})
can expressed as
$$
\kappa_i(A,b)=\frac{1}{\sigma_b}\left(
\frac{\neuc{C_i}^2}{\sigma_b^2}
\frac{\neuc{r}^2}{\alpha^2}
+c_{ii}(\frac{\neuc{x}^2}{\alpha^2}
+\frac{1}{\beta^2})\right)^{\frac{1}{2}}.
$$
The quantity $\sigma_b^2$ will often be estimated by $\frac{1}{m-n}\neuc{r}^2$
in which case the expression can be simplified
\begin{equation} \label{form:cii2}
\kappa_i(A,b)=\frac{1}{\sigma_b}\left(
\frac{1}{\alpha^2} \frac{\neuc{C_i}^2}{(m-n)}
+ \frac{c_{ii}\neuc{x}^2}{\alpha^2}
+ \frac{c_{ii}}{\beta^2})\right)^{\frac{1}{2}}.
\end{equation}
From Equation~(\ref{eq:normwise_formula}), we obtain
$$
\kappa_{LS}(A,b)=
   \frac{\neuc{C}^{1/2}}{\sigma_b}
   \left(
      \frac{
         \neuc{C} \neuc{r}^2 
      }
      { \alpha^2 \sigma_b^2 }
	  + \frac{\neuc{x}^2}{\alpha^2}
      + \frac{1}{\beta^2}
   \right)^{\frac{1}{2}}.
$$
The quantity $\sigma_b^2$ will often be estimated by $\frac{1}{m-n}\neuc{r}^2$
in which case the expression can be simplified
$$
\kappa_{LS}(A,b)=
   \frac{\neuc{C}^{1/2}}{\sigma_b}
   \left(
      \frac{ 1 } { \alpha^2 }
      \frac{ \neuc{C} } {  (m-n) }
	  + \frac{\neuc{x}^2}{\alpha^2}
      + \frac{1}{\beta^2}
   \right)^{\frac{1}{2}}.
$$
\section{Computation with LAPACK}\label{sec:lapack}
Section~\ref{sec:theoback} provides us with formulas to
compute the condition numbers $\kappa_i$ and $\kappa_{LS}$.  As explained in
Section~\ref{sec:link}, those quantities are intimately interrelated with the
entries of the variance-covariance matrix.  The goal of this section is to
present practical methods and codes to compute those quantities
efficiently with LAPACK.
The assumption made is that the LLSP has already been
solved with either the normal equations method or a QR factorization approach. Therefore
the solution vector $\hat x$, the norm of the residual $\| \hat r \|_2$, and
the R-factor $R$ of the QR factorization of $A$ are readily available (we recall
that the Cholesky factor of the normal equations is the R-factor of the QR
factorization up to some signs). In the example codes, we have used the
LAPACK routine DGELS
that solves the LLSP using QR factorization of A. Note that 
it is possible to have a more accurate solution using extra-precise 
iterative refinement~\cite{IRLLS.07}.
\subsection{Variance-covariance computation}
We will use the fact that $\frac{1}{m-n}\neuc{b-A\hat{x}}^2$ is an unbiased estimate
of $\sigma_b^2$. We wish to compute 
the following quantities related to the variance-covariance matrix $C$
\begin{itemize}
\item the $i$-th column $ C_i  =\sigma_b^2 e_i^T ( A^T A )^{-1}  $
\item the $i$-th diagonal element $c_{ii} =\sigma_b^2 \| e_i ^T A^\dagger \|_2^2 $
\item the whole matrix $ C $
\end{itemize}
We note that the quantities $ C_i $, $c_{ii}$, and $C$ are of interest
for statisticians. The NAG routine F04YAF~\cite{nag} is indeed an example
of tool to compute these three quantities.\\
For the two first quantities of interest, we note that
$$
\neuc{e_i^T A^{\dagger}}^2=\neuc{R^{-T}e_i}^2
~{\rm and}~
\neuc{e_i^T(A^TA)^{-1}}=\neuc{R^{-1}(R^{-T}e_i)}.
$$
\subsubsection{Computation of the $i$-th column $ C_i $}
$C_i$ can be computed with two $n$--by--$n$ triangular solves
\begin{equation}\label{form:trsys}
R^Ty=e_i~{\rm and}~ Rz=y.
\end{equation}
The $i$-th column of $C$ can be computed by
the following code fragment.\\
\\
{\bf
Code 1:\\
CALL DGELS(~'N',~M,~N,~1,~A,~LDA,~B,~LDB,~WORK,~LWORK,~INFO~)\\
RESNORM = DNRM2(~(M-N),~B(N+1),~1)\\
SIGMA2 = RESNORM**2/DBLE(M-N)\\
E(1:N) = 0.D0\\
E(I) = 1.D0\\
CALL DTRSV(~'U',~'T',~'N',~N-I+1,~A(I,I),~LDA,~E(I),~1)\\
CALL DTRSV(~'U',~'N',~'N',~N,~A,~LDA,~E,~1)\\
CALL DSCAL(~N,~SIGMA2,~E,~1)\\
}
\\
This requires about $2n^2$ flops (in addition to the cost of solving
the linear least squares problem using DGELS).\\
$c_{ii}$ can be computed by one $n$--by--$n$ triangular solve and taking the
square of the norm of the solution which involves about $(n-i+1)^2$ flops. It is
important to note that the larger $i$, the less expensive to obtain $c_{ii}$. In
particular if $i=n$ then only one operation is needed: $c_{nn} = R_{nn}^{-2}$.
This suggests that a correct ordering of the variables can save some
computation.\\
\\
\subsubsection{Computation of the $i$-th diagonal element $ c_{ii} $}
From $c_{ii} = \sigma_b^2 \neuc{e_i^TR^{-1}}^2$, it comes that each $ c_{ii} $
corresponds to the $i$-th row of $R^{-1}$.
Then the diagonal elements of $C$ can be computed by the
following code fragment.\\
\\
{\bf
Code 2:\\
CALL DGELS(~'N',~M,~N,~1,~A,~LDA,~B,~LDB,~WORK,~LWORK,~INFO~)\\
RESNORM~=~DNRM2((M-N), B(N+1), 1)\\
SIGMA2~=~RESNORM**2/DBLE(M-N)\\
CALL DTRTRI(~'U',~'N',~N,~A,~LDA,~INFO)\\
DO~I=1,N \\
\hspace*{0.5cm} CDIAG(I)~=~DNRM2(~N-I+1,~A(I,I),~LDA)\\
\hspace*{0.5cm} CDIAG(I)~=~SIGMA2~*~CDIAG(I)**2\\
END DO\\
}
\\
This requires about $n^3/3$ flops (plus the cost of DGELS).\\
\\
\subsubsection{Computation of the whole matrix $C$}
In order to compute explicity all the coefficients of the matrix $C$, one can
use the routine DPOTRI which computes the inverse of a matrix from its Cholesky
factorization. First the routine computes the inverse of $R$ using DTRTRI
and then performs the triangular matrix-matrix multiply $R^{-1}R^{-T}$ by
DLAUUM. This requires about $2n^3/3$ flops.
We can also compute the variance-covariance matrix without inverting $R$
using for instance the algorithm given in~\cite[p. 119]{BJORCK} but the
computational cost remains $2n^3/3$ (plus the cost of DGELS).\\
\\
We can obtain the upper triangular part of $C$
by the following code fragment.\\
\\
{\bf
Code 3:\\
CALL DGELS(~'N',~M,~N,~1,~A,~LDA,~B,~LDB,~WORK,~LWORK,~INFO~)\\
RESNORM~=~DNRM2((M-N), B(N+1), 1)\\
SIGMA2~=~RESNORM**2/DBLE(M-N)\\
CALL DPOTRI(~'U',~N,~A,~LDA,~INFO)\\
CALL DLASCL(~'U',~0,~0,~N,~N,~1.D0,~SIGMA2,~N,~N,~A,~LDA,~INFO)\\
}
\\
\subsection{Condition numbers computation}
For computing $\kappa_i(A,b)$, we need to compute both
the $i$-th diagonal element and the norm of the $i$-th column
of the variance-covariance matrix and we cannot use direcly Code 1 but the
following code fragment\\
\\
{\bf
Code 4:\\
ALPHA2 = ALPHA**2\\
BETA2 = BETA**2\\
CALL DGELS(~'N',~M,~N,~1,~A,~LDA,~B,~LDB,~WORK,~LWORK,~INFO~)\\
XNORM = DNRM2(N,~B(1),~1)\\
RESNORM = DNRM2((M-N), B(N+1), 1)\\
CALL DTRSV(~'U',~'T',~'N',~N-I+1,~A(I,I),~LDA,~E(I),~1~)\\
ENORM = DNRM2(N, E, 1)\\
K = (ENORM**2)*(XNORM**2/ALPHA2+1.d0/BETA2)\\
CALL DTRSV(~'U',~'N',~'N',~N,~A,~LDA,~E,~1~)\\
ENORM = DNRM2(N, E, 1)\\
K = SQRT((ENORM*RESNORM)**2/ALPHA2 + K)\\
}
\\
For computing all the $\kappa_i(A,b)$,
we need to compute the columns $C_i$ and the diagonal elements $c_{ii}$ 
using Formula~(\ref{form:cii2}) and then we have to compute 
the whole variance-covariance matrix. This can be 
performed by a slight modification of Code 3.\\
When only $b$ is perturbed, then we have to invert $R$ and we can use 
a modification of Code 2 (see numerical example in Section~\ref{sec:cnes}).\\
\\
For estimating $\kappa_{LS}(A,b)$, we need to have an estimate of
$\neuc{R^{-1}}$.
The computation of $\neuc{R^{-1}}$ requires to compute the minimum singular value of the 
matrix $A$ (or $R$). 
One way is to compute the full SVD of $A$ (or $R$) which requires $\bigo (n^3)$ flops.
As an alternative,
$\neuc{R^{-1}}$ can be estimated for instance by considering
other matrix norms through the following inequalities
\begin{eqnarray*}
\frac{1}{\sqrt{n}}\nfro{R^{-1}} & \leq \neuc{R^{-1}} \leq & \nfro{R^{-1}}, \\
\frac{1}{\sqrt{n}}\|R^{-1}\|_{\infty} & \leq \neuc{R^{-1}} \leq & \sqrt{n} \|R^{-1}\|_{\infty},\\
\frac{1}{\sqrt{n}}\|R^{-1}\|_1 & \leq \neuc{R^{-1}} \leq & \sqrt{n} \|R^{-1}\|_1.\\
\end{eqnarray*}
$\|R^{-1}\|_1$ or $\|R^{-1}\|_{\infty}$ can be estimated using
Higham modification~\cite[p. 293]{HIGHAM} of Hager's~\cite{HAGER}
method as it is implemented in LAPACK~\cite{LAPACK} DTRCON routine
(see Code 5). The cost is $\bigo( n^2 )$.\\
\\
{\bf 
Code 5:\\
CALL DTRCON(~'I',~'U',~'N',~N,~A,~LDA,~RCOND,~WORK,~IWORK,~INFO)\\ 
RNORM~=~DLANTR(~'I',~'U',~'N',~N,~N,~A,~LDA,~WORK)\\
RINVNORM~=~(1.D0/RNORM)/RCOND\\
}
\\
We can also evaluate $\neuc{R^{-1}}$ by considering
$\nfro{R^{-1}}$ since we have
\begin{eqnarray*}
\nfro{R^{-1}}^2 & = & \nfro{R^{-T}}^2\\
& = & {\rm tr}(R^{-1}R^{-T})\\
& = & \frac{1}{\sigma_b^2} {\rm tr}(C),\\
\end{eqnarray*}
where tr($C$) denotes the trace of the matrix $C$, i.e
$\sum_{i=1}^{n} c_{ii}$.
Hence the condition number of the
least-squares solution can be approximated by 
\begin{equation} \label{form:trace}
\kappa_{LS}(A,b) \simeq
\left( \frac{{\rm tr}(C)}{\sigma_b^2} \left(\frac{{\rm tr}(C)
\neuc{r}^2+\sigma_b^2 \neuc{x}^2}
{\sigma_b^2 \alpha^2}+\frac{1}{\beta^2}\right) \right)^{\frac{1}{2}}.
\end{equation}
Then we can estimate $\kappa_{LS}(A,b)$ by computing and summing
the diagonal elements of $C$ using Code 2.\\

When only $b$ is perturbed ($\alpha = +\infty~{\rm and}~\beta=1$),
then we get
$$\kappa_{LS}(b) \simeq \frac{\sqrt{{\rm tr}(C)}}{\sigma_b}.$$
This result relates to~\cite[p. 167]{FAREBROTHER} where
${\rm tr}(C)$ measures the squared effect on the
LLSP solution $x$ to small changes in $b$.\\
\\
We give in Table~\ref{TabCompar} the LAPACK routines used for
computing the condition numbers of an LLSP solution or
its components and the corresponding number of floating-point operations
per second. Since the LAPACK routines involved in the covariance and/or
LLSP condition numbers have their equivalent in the parallel
library ScaLAPACK~\cite{SCALAPACK}, then this table is also available 
when using ScaLAPACK. 
This enables us to easily compute these quantities for larger LLSP.

\begin{table}[hbtp!]
\centering
\caption{Computation of least squares conditioning with (Sca)LAPACK}
\vspace{0.4cm}
\begin{tabular}{|c|c|c|c|}
\hline
condition number&linear algebra operation&LAPACK routines&flops count\\
&&&\\
\hline
$\kappa_i(A,b)$&$R^Ty=e_i~{\rm and}~ Rz=y$& 2 calls to (P)DTRSV&$2n^2$\\
&&&\\
\hline
all $\kappa_i(A,b),~i=1,n$&$RY=I~{\rm and~compute}~YY^T$&(P)DPOTRI&$2n^3/3$\\
&&&\\
\hline
all $\kappa_i(b),~i=1,n$&invert $R$&(P)DTRTRI&$n^3/3$\\
&&&\\
\hline
$\kappa_{LS}(A,b)$&estimate $\|R^{-1}\|_{1~{\rm or}~\infty}$&(P)DTRCON&${\cal O}(n^2)$\\
&compute $\nfro{R^{-1}}$&(P)DTRTRI&$n^3/3$\\
\hline
\end{tabular}
\label{TabCompar}
\end{table}
\begin{Remark}
{\em
The cost for computing all the $\kappa_i(A,b)$ or estimating
$\kappa_{LS}(A,b)$ is always ${\cal O}(n^3)$.
This seems affordable when we compare it to the cost of the
least squares solution using Householder QR factorization ($2mn^2-2n^3/3$)
or the normal equations ($mn^2+n^3/3$) because we have in general $m \gg n$.\\
}
\end{Remark}
\section{Numerical experiments}\label{sec:numerics}
\subsection{Laplace's computation of the mass of Jupiter and assessment of
the validity of its results}
In~\cite{LAPLACE}, Laplace computes the mass of Jupiter, Saturn and Uranus
and provides the variances associated with those variables in order to
assess the quality of the results. The data comes from the French
astronomer Bouvart in the form of the normal equations given in
Equation~(\ref{eq:laplace}).
\begin{equation}\label{eq:laplace}
\begin{array}{rcl}
    795938 z_0 -  12729398 z_1 +   6788.2 z_2 -  1959.0 z_3 + 696.13 z_4 +
2602 z_5 & = &    7212.600 \\
 -12729398 z_0 + 424865729 z_1 - 153106.5 z_2 - 39749.1 z_3 -   5459 z_4 +
5722 z_5 & = & -738297.800 \\
    6788.2 z_0 -  153106.5 z_1 +  71.8720 z_2 -  3.2252 z_3 + 1.2484 z_4 +
1.3371 z_5 & = &     237.782 \\
   -1959.0 z_0 -   39749.1 z_1 -   3.2252 z_2 + 57.1911 z_3 + 3.6213 z_4 +
1.1128 z_5 & = &     -40.335 \\
    696.13 z_0 -      5459 z_1 +   1.2484 z_2 +  3.6213 z_3 + 21.543 z_4 +
46.310 z_5 & = &    -343.455 \\
      2602 z_0 +      5722 z_1 +   1.3371 z_2 +  1.1128 z_3 + 46.310 z_4 +
129 z_5 & = &   -1002.900 \\
\end{array}
\end{equation}
For computing the mass of Jupiter, we know that Bouvart performed $ m =
129 $ observations and there are $n=6$ variables in the system. The
residual of the solution $\| b - A\hat x \|_2^2$ is also given by Bouvart
and is $31096$.  On the $6$ unknowns, Laplace only seeks one, the second
variable $z_1$.  The mass of Jupiter in term of the mass of the Sun is
given by $z_1$ and the formula:
$$ \textmd{mass of Jupiter} = \frac{1+z_1}{1067.09}.$$
It turns out that the first variable $z_0$
represents the mass of Uranus through the formula
$$ \textmd{mass of Uranus} = \frac{1+z_0}{19504}.$$
If we solve the system~(\ref{eq:laplace}), we obtain the solution vector\\
\begin{center}
\begin{minipage}{10cm}
Solution vector\\
0.08954  -0.00304 -11.53658  -0.51492   5.19460 -11.18638
\end{minipage}
\end{center}
From $z_1$, we can compute the mass of Jupiter as a fraction of the mass of the Sun
and we obtain $1070$.
This value is indeed accurate since the correct value according to NASA is
$1048$. From $z_0$, we can compute the mass of Uranus as a fraction of the mass of the Sun and we obtain
$17918$. This value is inaccurate since the correct value according to NASA is
$22992$.\\
Laplace has computed the variance of $z_0$ and $z_1$ to assess the fact that
$z_1$ was probably correct and $z_0$ probably inaccurate. To compute those
variances, Laplace first performed a Cholesky factorization from right to left of
the system~(\ref{eq:laplace}), then, since the variables were correctly ordered
the number of operations involved in the computation of the variances of $z_0$ and $z_1$
were minimized. The variance-covariance matrix for Laplace's system is:
$$
\left(
\begin{array}{cccccc}
0.005245 & -0.000004 & -0.499200 &  0.137212 &   0.235241 &  -0.186069 \\
\cdot    &  0.000004 &  0.009873 &  0.003302 &   0.002779 &  -0.001235 \\
\cdot    & \cdot     & 71.466023 & -5.441882 & -16.672689 &  14.922752 \\
\cdot    & \cdot     & \cdot     & 10.860492 &   5.418506 &  -4.896579 \\
\cdot    & \cdot     & \cdot     & \cdot     &  66.088476 & -28.467391 \\
\cdot    & \cdot     & \cdot     & \cdot     &  \cdot     &  15.874809 \\
\end{array}
\right)
$$

Our computation gives us that the variance for the mass of Jupiter is
$4.383233\cdot10^{-6}$. For reference, Laplace in 1820 computed
$4.383209\cdot10^{-6}$. (We deduce the variance from Laplace's value 5.0778624.
To get what we now call the variance, one needs to compute the quantity:
$ 1/(2*10**5.0778624)*m/(m-n)$.)

From the variance-covariance matrix, one can assess that the computation
of the mass of Jupiter (second variable) is extremely reliable while the
computation of the mass of Uranus (first variable) is not.
For more details, we recommend to read \cite{langoureview2007}.

\subsection{Gravity field computation}\label{sec:cnes}
\begin{figure}[!ht]
      \begin{center}
      {\epsfig{file=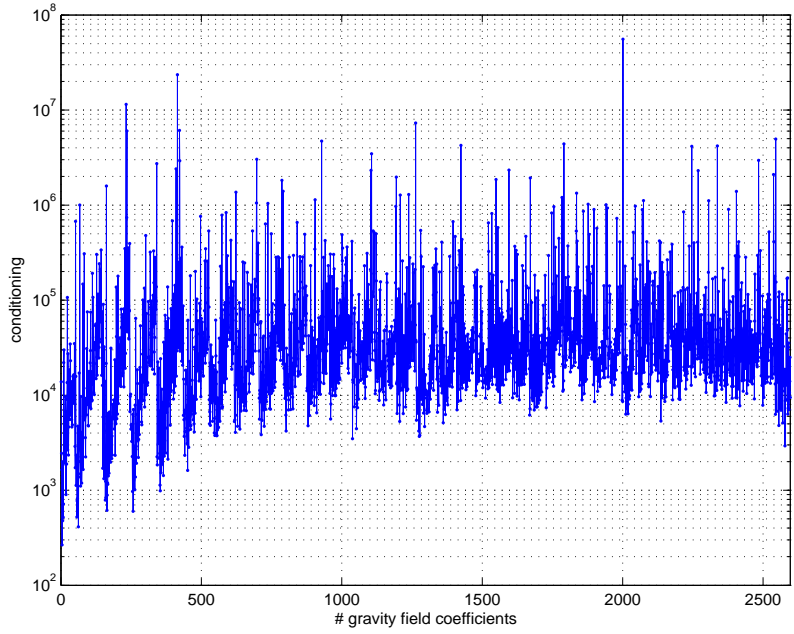,width=0.7 \textwidth}} \\
      \caption{\label{plotcond}
        Amplitude of the relative condition numbers for the gravity field coefficients.
      }
     \end{center}
\end{figure}

\begin{figure}[!ht]
      \begin{center}
      {\epsfig{file=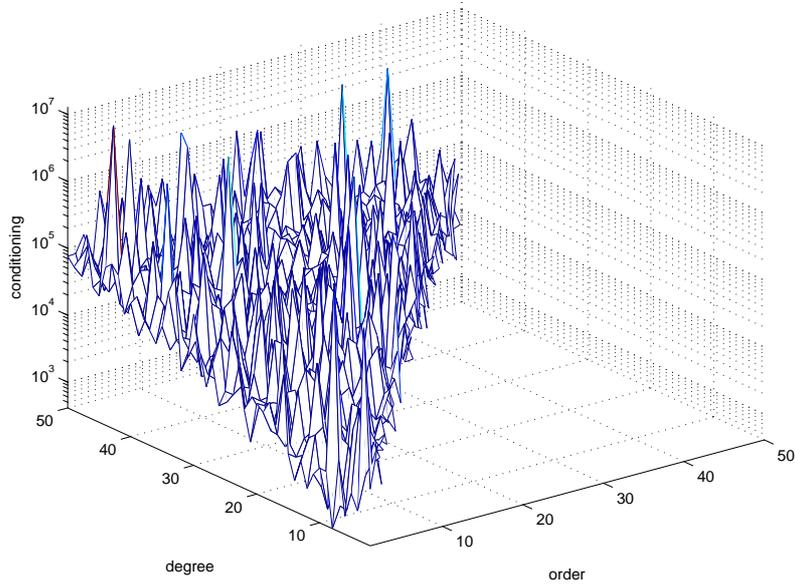,width=0.7 \textwidth}} \\
      \caption{\label{plotcos}
        Conditioning of spherical harmonic coefficients $\overline{C}_{\ell m}~(2 \leq \ell \leq 50~,~1 \leq m\leq 50)$.
      }
     \end{center}
\end{figure}

\begin{figure}[!ht]
      \begin{center}
      {\epsfig{file=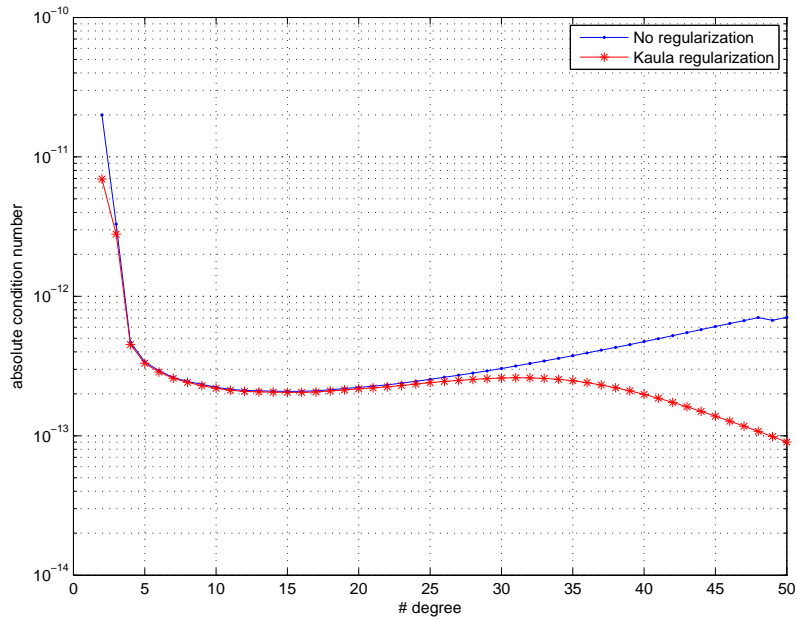,width=0.7 \textwidth}} \\
      \caption{\label{zonaux}
Effect of regularization on zonal coefficients $\overline{C}_{{\ell} 0}~(2 \leq {\ell} \leq 50)$        
      }
     \end{center}
\end{figure}

A classical example of parameter estimation problem is the computation of the 
Earth's gravity field coefficients. More specifically, we estimate the 
parameters of the gravitational potential that can be expressed in spherical 
coordinates $(r,\theta,\lambda)$ by~\cite{BALMINO}
\begin{equation} \label{potential}
V(r,\theta,\lambda)=\frac{GM}{R}\sum_{\ell=0}^{\ell_{max}}\left(\frac{R}{r}\right)^{\ell+1}
\sum_{m=0}^{\ell}\overline{P}_{\ell m}(\cos{\theta})\left[\overline{C}_{\ell m}
\cos{m\lambda}+\overline{S}_{\ell m}\sin{m\lambda}\right]
\end{equation}
where $G$ is the gravitational constant, $M$ is the Earth's mass,
$R$ is the Earth's reference radius,
the $\overline{P}_{\ell m}$ represent the fully normalized Legendre
functions of degree $\ell$ and
order $m$ and $\overline{C}_{\ell m}$,$\overline{S}_{\ell m}$ 
are the corresponding normalized harmonic coefficients.
The objective here is to compute the harmonic coefficients $\overline{C}_{\ell m}$ and
$\overline{S}_{\ell m}$ the most accurately as possible.
The number of unknown parameters is expressed by $n=(\ell_{max}+1)^2.$
These coefficients are computed by solving a linear least squares 
problem that may involve millions of observations and tens of thousands of variables.
More details about the physical problem and the resolution methods can be found in
\cite{PHD.MB}.
The data used in the following experiments were provided by 
CNES\footnote{Centre National d'Etudes Spatiales, Toulouse, France}
and they correspond to 10 days of observations
using GRACE\footnote{Gravity Recovery and
Climate Experiment, NASA, launched March 2002} measurements 
(about $166,000$ observations). 
We compute the spherical harmonic coefficients 
$\overline{C}_{\ell m}$ and $\overline{S}_{\ell m}$ up to a degree 
$\ell_{max}=50$; except the coefficients 
$\overline{C}_{11}, \overline{S}_{11}, \overline{C}_{00}, \overline{C}_{10}$ 
that are a priori known. 
Then we have $n=2,597$ unknowns in the corresponding least squares problems
(note that the GRACE satellite enables us to compute 
a gravity field model up to degree 150).
The problem is solved using the normal equations method and we have the Cholesky 
decomposition $A^TA=U^TU$.\\ 
We compute the relative condition numbers of each coefficient $x_i$ 
using the formula 
$$\kappa^{(rel)}_i(b)=\neuc{e_i^TU^{-1}} \neuc{b}/|x_i|,$$ 
and the following code fragment, derived from Code 2, in which
the array $D$ contains the normal equations $A^TA$ and the vector $X$ 
contains the right-hand side $A^Tb$.\\
\\
{\bf CALL DPOSV(~'U',~N,~1,~D,~LDD,~X,~LDX,~INFO)\\
CALL DTRTRI(~'U',~'N',~N,~D,~LDD,~INFO)\\
DO~I=1,N \\
\hspace*{0.5cm} KAPPA(I)~=~DNRM2(~N-I+1,~D(I,I),~LDD)~*~BNORM/ABS(X(I)) \\
END DO}\\
\\
Figure~\ref{plotcond} represents the relative condition numbers of all the $n$ 
coefficients. We observe the disparity between the condition numbers   
(between $10^2$ and $10^8$).
To be able to give a physical interpretation, we need first to sort the coefficients 
by degrees and orders as given in the development of $V(r,\theta,\lambda)$  
in Expression~(\ref{potential}).\\
In Figure~\ref{plotcos}, we plot the coefficients $\overline{C}_{\ell m}$ as 
a function of the degrees and orders (the curve with the $\overline{S}_{\ell m}$ 
is similar). We notice that for a given order, the condition number increases 
with the degree and that, for a given degree, the variation of the sensitivity 
with the order is less significant.\\
We can also study the effect of regularization on the conditioning.
The physicists use in general a Kaula~\cite{KAULA} regularization technique 
that consists of adding to $A^TA$ a diagonal matrix 
$D=diag(0,\cdots,0,\delta,\cdots,\delta)$ 
where $\delta$ is a constant that is proportional to 
$\frac{10^{-5}}{\ell_{max}^2}$ and 
the nonzero terms in $D$ correspond to the variables 
that need to be regularized. An example of the effect of Kaula regularization 
is shown in Figure~\ref{zonaux} where we consider
the coefficients of order $0$ also called zonal coefficients.
We compute here the absolute condition numbers of these coefficients 
using the formula $\kappa_i(b)=\neuc{e_i^TU^{-1}}$.
Note that the $\kappa_i(b)$ are much lower that 1. This is not surprising
because typically in our application $\neuc{b} \sim 10^5$/ and $|x_i| \sim 10^{-12}$ 
which would make the associated relative condition numbers greater than 1.
We observe that the regularization 
is effective on coefficients of highest degree that are in general more sensitive 
to perturbations. 

\newpage

\section{Conclusion}\label{sec:concl}
To assess the accuracy of a linear least squares solution, the practitioner of
numerical linear algebra uses generally quantities like condition
numbers or backward errors when the statistician is more interested in
covariance analysis.
In this paper we proposed quantities that talk to both communities
and that can assess the quality of the solution of a least squares problem
or one of its component. We provided pratical ways to compute these quantities
using (Sca)LAPACK and we experimented these computations
on pratical examples including a real physical application in the area of space geodesy.
\bibliographystyle{siam}
\bibliography{biblio}
\end{document}